\documentclass[a4paper,12pt]{amsart}
\usepackage{hyperref,fancyhdr,mathrsfs,amsmath,amscd,amsthm,amsfonts,latexsym,amssymb,stmaryrd}

\DeclareMathOperator{\trace}{trace}
\DeclareMathOperator{\grad}{grad}

\DeclareMathOperator{\dif}{d}

\DeclareMathOperator{\Lie}{\mathcal{L}}

\renewcommand{\H}{\mathscr{H}}
\newcommand{\V}{\mathscr{V}}
\newcommand{\F}{\mathscr{F}}

\def \a{\alpha}

\def \G{\Gamma}
\def \l{\lambda}

\def \O{\Omega}
\def \phi{\varphi}
\def \Phi{\varPhi}
\def \p{\pi}

\def \s{\sigma}

\def \R{\mathbb{R}}

\def \C{\mathbb{C}\,}

\def\widecheckg{g^{\hspace*{-2.5pt}\vbox to 5pt{\hbox to
0pt{\LARGE$\check{}$}}}\hspace*{2pt}}

\def\widecheckl{\lambda^{\hspace*{-3.5pt}\vbox to 8pt{\hbox to
0pt{\LARGE$\check{}$}}}\hspace*{2pt}}

\begin{document}

\title{Harmonic morphisms with one-dimensional fibres on conformally-flat Riemannian manifolds}
\author{Radu~Pantilie}
\thanks{The author gratefully acknowledges that this work was partially supported by a
CEx Grant no.\ 2-CEx 06-11-22/25.07.2006.}
\email{\href{mailto:Radu.Pantilie@imar.ro}{Radu.Pantilie@imar.ro}}
\address{R.~Pantilie, Institutul de Matematic\u a ``Simion Stoilow'' al Academiei Rom\^ane,
C.P. 1-764, 014700, Bucure\c sti, Rom\^ania}

\newtheorem{thm}{Theorem}[section]
\newtheorem{lem}[thm]{Lemma}
\newtheorem{cor}[thm]{Corollary}
\newtheorem{prop}[thm]{Proposition}

\theoremstyle{definition}

\newtheorem{defn}[thm]{Definition}
\newtheorem{rem}[thm]{Remark}
\newtheorem{exm}[thm]{Example}

\numberwithin{equation}{section}

\maketitle
\thispagestyle{empty}
\vspace{-0.5cm}
\section*{Abstract}
\begin{quote}
{\footnotesize  We classify the harmonic morphisms with one-dimensional fibres
(1) from real-analytic conformally-flat Riemannian manifolds of dimension at least four,
and (2) between conformally-flat Riemannian manifolds of dimensions at least
three.}
\end{quote}

\section*{Introduction}

\indent
Harmonic morphisms between Riemannian manifolds are maps which pull-back (local)
harmonic functions to harmonic functions. By a basic result, a map is a harmonic morphisms
if and only if it is a harmonic map which is horizontally weakly conformal (see \cite{BaiWoo2}\,).\\
\indent
There are, now, several classification results for harmonic morphisms with one-dimensional fibres.
In \cite{Bry}\,, it was proved that there are precisely two types of such harmonic morphisms
from Riemannian manifolds, with constant curvature, of dimension at least four.
This result was generalized, in \cite{PanWoo-d}\,, to Einstein manifolds of dimension
at least five; in dimension four, the situation is different, there appears
a third type of harmonic morphism \cite{Pan-4to3} (see \cite{PanWoo-exm}\,).
Also, in \cite{PanWoo-sd}\,, are classified the `twistorial' harmonic morphisms
with one-dimensional fibres from self-dual four-manifolds.\\
\indent
In this paper, we classify the harmonic morphisms with
one-dimensional fibres from conformally-flat Riemannian manifolds
of dimension at least four. We prove  that there are just two types
of such harmonic morphisms, one of which (the `Killing type'), also,
appears in the above mentioned results, whilst the second type is an extension
of the `warped product type', involved in \cite{Bry}\,, \cite{Pan-4to3}\,,
\cite{PanWoo-d} and \cite{PanWoo-sd}\,.\\
\indent
The main result is given in Section \ref{section:main} (Theorem \ref{thm:cf1d}\,) after
a brief review of harmonic morphisms, and conformally-flat Riemannian manifolds, given
in Sections \ref{section:harmorphs} and \ref{section:conf-flat}\,, respectively.\\
\indent
In Section \ref{section:main}\,, we also classify the harmonic morphisms with one-dimensional
fibres between conformally-flat Riemannian manifolds (Corollary \ref{cor:cf1d}\,).
It follows that \emph{the Hopf polynomial map $\R^4\to\R^3$\,,
$(z_1,z_2)\mapsto(|z_1|^2-|z_2|^2,2z_1\overline{z_2})$\,, is,
up to local conformal diffeomorphisms, the only harmonic morphism with one-dimensional fibres
and nonintegrable horizontal distribution between conformally-flat Riemannian manifolds,
of dimensions at least three} (Corollary \ref{cor:cf1dHopf}\,).\\
\indent
I am grateful to John~C.~Wood for useful comments.

\section{Harmonic morphisms with one-dimensional fibres} \label{section:harmorphs}

\indent
In this section we recall a few facts on harmonic morphisms with one-dimensional fibres.\\
\indent
Unless otherwise stated, all the manifolds are assumed to be connected and smooth
and all the maps are assumed to be smooth.

\begin{defn}
A \emph{harmonic morphism (between Riemannian manifolds)} is a map $\phi:(M^m,g)\to(N^n,h)$
such that if $U$ is an open set of $N$, with $\phi^{-1}(U)\neq\emptyset$\,, and
$f$ is a harmonic function on $(U,h|_U)$ then $f\circ\phi$ is a harmonic function
on $(\phi^{-1}(U),g|_{\phi^{-1}(U)})$\,.
\end{defn}

\begin{defn}
A map between Riemannian manifolds $\phi:(M^m,g)\to(N^n,h)$ is \emph{horizontally
weakly conformal} if, for any $x\in M$, either $\dif\!\phi_x=0$ or,
for any $X,Y\in({\rm ker}\dif\!\phi_x)^{\perp}$, we have
$h(\dif\!\phi(X),\dif\!\phi(Y))=\l(x)^2g(X,Y)$ for some positive number $\l(x)$\,.\\
\indent
The function $\l$\,, extended to be zero over the set of points $x\in M$ where $\dif\!\phi_x=0$\,,
is the \emph{dilation} of $\phi$\,. (Note that the dilation $\l$ is continuous on $M^m$ whilst
the \emph{square dilation} $\l^2$ is smooth on $M^m$.)
\end{defn}

\indent
Let $\phi:(M^m,g)\to(N^n,h)$ be a horizontally conformal submersion; denote by $\l$ its dilation.
Then $\l=1$ if and only if $\phi$ is a \emph{Riemannian submersion}.\\

\indent
If $m=n$ then a surjective map $\phi:(M^m,g)\to(N^m,h)$ is a horizontally conformal submersion
if and only if it is a \emph{conformal diffeomorphism}. The dilation of a conformal diffeomorphism $\phi$
is called the \emph{conformality factor} of $\phi$\,.\\

\indent
The study of harmonic morphisms is based on the following result of B.~Fuglede and T.~Ishihara
(see \cite{BaiWoo2}\,).

\begin{thm} \label{thm:FugIsh}
A map is a harmonic morphism if and only if it is a harmonic map which is horizontally
weakly conformal.
\end{thm}

\indent
As usual, if $\phi:(M^m,g)\to N^n$ is a submersion we denote by $\V={\rm ker}\dif\!\phi$
the \emph{vertical distribution} and by $\H=\V^{\perp}$ the \emph{horizontal distribution}.\\
\indent
For horizontally conformal submersions with one-dimensional fibres the condition
of harmonicity can be expressed as follows.

\begin{prop}[\,\cite{Bry}\,, see \cite{Pan}\,, \cite{BaiWoo2}\,] \label{prop:Bryant}
Let $\phi:(M^{n+1},g)\to(N^n,h)$ be a horizontally conformal submersion with one-dimensional
fibres; $n\geq3$\,. Let $\l$ be the dilation of $\phi$ and let $V$ be the vertical vector field
(well-defined up to sign) such that $g(V,V)=\l^{2n-4}$.\\
\indent
The following assertions are equivalent.\\
\indent
\quad{\rm (i)} $\phi:(M^{n+1},g)\to(N^n,h)$ is a harmonic morphism.\\
\indent
\quad{\rm (ii)} $[V,X]=0$ for any basic (horizontal) vector field $X$\,.\\
\indent
Furthermore, if {\rm (i)} or {\rm (ii)} holds then $\O=\dif\!\theta$ is basic,
where $\theta$ is the vertical dual of $V$, characterised by $\theta(V)=1$ and
$\theta|_{\H}=0$\,.
\end{prop}

\indent
Let $\phi:(M^{n+1})\to(N^n,h)$ be a harmonic morphism with one-dimensional fibres.
With the same notations as in Proposition \ref{prop:Bryant}\,, the vector field $V$
is called the \emph{fundamental vector field}. It is easy to prove that $\O=0$
if and only if $\H$ is integrable; also, we have $g=\l^{-2}\phi^*(h)+\l^{2n-4}\theta^2$.
It follows that $V$ is a Killing vector field if and only if $\V$ is a Riemannian foliation;
equivalently, $V(\l)=0$ (\,\cite{Bry}\,, see \cite{Pan}\,, \cite{BaiWoo2}\,).
If $V$ is Killing then $\phi$ is called of \emph{Killing type}.

\begin{lem}[\,\cite{Pan-thesis}\,, cf.\ \cite{BaiWoo2}\,] \label{lem:curv}
Let $\phi:(M^{n+1},g)\to(N^n,h)$ be a harmonic morphism. Let $\l$ be its dilation
and let $V$ be the fundamental vector field of $\phi$\,; we shall denote by
$\s=\log\l$\,.\\
\indent
Then we have the following relations for the curvature tensors
$R^M$ and $R^N$ of $(M^{n+1},g)$ and $(N^n,h)$\,, respectively:
\begin{equation} \label{e:xvyv} \begin{split}
R^M(X,V,Y,V)=&-\frac{1}{2}(n-2)e^{(2n-4)\s}(\Lie_{\H(\grad_{h}\s)}h)(X,Y)\\
&-(n-2)e^{(2n-4)\s}\{nX(\s)Y(\s)-|\H(\grad_{h}\s)|_{h}^{2}\;h(X,Y)\}\\
&+e^{-2\s}\{V(V(\s))-(n-1)V(\s)^{2}\}h(X,Y)\\
&+\frac{1}{4}\,e^{(4n-6)\s}\;h(i_{X}\O,i_{Y}\O)\:,
\end{split} \end{equation}
\begin{equation} \label{e:xyzv} \begin{split}
R^M(X,&Y,Z,V)=-\frac{1}{2}e^{(2n-4)\s}(^{h}\nabla\O)(X,Y,Z)\\
&+\frac{1}{2}(n-1)e^{(2n-4)\s}\{X(\s)\O(Y,Z)+Y(\s)\O(Z,X)-2Z(\s)\O(X,Y)\}\\
&-e^{-2\s}\{X(V(\s))-(n-2)X(\s)V(\s)\}h(Y,Z)\\
&+e^{-2\s}\{Y(V(\s))-(n-2)Y(\s)V(\s)\}h(X,Z)\\
&+\frac{1}{2}e^{(2n-4)\s}\{\O(X,\grad_{h}\s)h(Y,Z)-\O(Y,\grad_{h}\s)h(X,Z)\}\:,
\end{split} \end{equation}
\begin{equation} \label{e:xyzh} \begin{split}
&R^M(X,Y,Z,H)=e^{-2\s}\phi^*(R^N)(X,Y,Z,H)\\
&-\frac14\,e^{(2n-4)\s}\{\O(H,X)\O(Y,Z)+\O(H,Y)\O(Z,X)-2\O(H,Z)\O(X,Y)\}\\
&-\frac12\,e^{-2\s}V(\s)\,\{-\O(Y,H)h(X,Z)+\O(X,H)h(Y,Z)-\O(X,Z)h(Y,H)+\O(Y,Z)h(X,H)\}\\
&-e^{-2\s}\{X(\s)H(\s)h(Y,Z)-X(\s)Z(\s)h(Y,H)-Y(\s)H(\s)h(X,Z)+Y(\s)Z(\s)h(X,H)\}\\
&+e^{-2\s}\{h(X,Z)h(^h\nabla_Y(\H(\grad_h\s)),H)-h(Y,Z)h(^h\nabla_X(\H(\grad_h\s)),H)\\
&+h(Y,H)h(^h\nabla_X(\H(\grad_h\s)),Z)-h(X,H)h(^h\nabla_Y(\H(\grad_h\s)),Z)\}\\
&-e^{-2\s}\{h(X,Z)h(Y,H)-h(X,H)h(Y,Z)\}\{e^{(-2n+2)\s}\,V(\s)^2+|\H(\grad_h\s|_h^2\}\:,
\end{split} \end{equation}
where $X,Y,Z,H$ are horizontal and $^{h}\nabla$ denotes the Levi-Civita connection of $(M,h)$.
\end{lem}

\begin{rem}
See \cite{BaiWoo2} and the references therein for more information on harmonic morphisms
between Riemannian manifolds and, in particular, for the notion of $p$-harmonic morphism.
Also, see \cite{LouPan,LouPan-II,Pan-Deva} for harmonic morphisms in the more general
setting of Weyl geometry.
\end{rem}

\section{Conformally-flat Riemannian manifolds} \label{section:conf-flat}

\indent
Firstly, we recall (see \cite{Laf-Weyl}\,) the definition of the Weyl tensor of a Riemannian
manifold.\\
\indent
Let $(M^m,g)$ be a Riemannian manifold. For $h$ and $k$ sections of $\odot^2(T^*M)$
(that is, $h$ and $k$ are symmetric covariant tensor fields
of degree two on $M^m$), we shall denote by $h\owedge k$ the section of
$\odot^2(\Lambda^2(T^*M))$ defined by
\begin{equation*}
\begin{split}
(h\owedge k)(T,X,Y,Z)=&h(T,Y)k(X,Z)+h(X,Z)k(T,Y)\\
&-h(T,Z)k(X,Y)-h(X,Y)k(T,Z)\;,
\end{split}
\end{equation*}
for any $T,X,Y,Z\in TM$.\\
\indent
If $S$ is a (1,3)-tensor field on $(M,g)$ then we shall denote by the same symbol $S$
the (0,4)-tensor field defined by $S(T,X,Y,Z)=-g(S(T,X,Y),Z)$\,, for any $T,X,Y,Z\in TM$.\\
\indent
The \emph{Weyl (curvature) tensor} of $(M^m,g)$ is the (1,3)-tensor field $W$ characterised
by the following two conditions:\\
\indent
\quad1) $\trace(X\mapsto W(X,Y)Z)=0$\,, for any $Y,Z\in TM$,\\
\indent
\quad2) $R=g\owedge r+W$ for some (necessarily unique) section $r$ of $\odot^2(T^*M)$\,,
where $R$ is the curvature tensor of $(M^m,g)$\,.\\
\indent
The Weyl tensor is conformally invariant; that is, if we denote by $W^g$ the Weyl tensor
of $(M^m,g)$ then $W^{\l^2g}=W^g$, for any positive function $\l$ on $M^m$.\\
\indent
The Riemannian manifold $(M^m,g)$ is called \emph{(locally) conformally-flat} if
for each point of $M^m$ there exists an open neighbourhood $U$ and a conformal diffeomorphism
$\phi$ from $U$ onto some open set of $\R^m$ (endowed with its canonical Riemannian metric);
the local coordinates on $U$ induced by $\phi$ are called \emph{flat}.\\
\indent
{}From Liouville's theorem on local conformal diffeomorphisms between Euclidean spaces (see \cite{BaiWoo2}\,)\,,
it follows easily that if $(M^m,g)$ is conformally-flat then $M^m$ is real-analytic in flat
local coordinates $(m\geq2)$\,.\\
\indent
The following theorem is due to H.~Weyl (see \cite{Laf-Weyl}\,).

\begin{thm} \label{thm:Weyl}
A Riemannian manifold, of dimension at least four, is conformally-flat if and only if its
Weyl tensor is zero.
\end{thm}
\noindent
(See \cite{Laf-Weyl} for the case when the dimension is less than four.)\\

\indent
We do not imagine that the following result is new.

\begin{prop} \label{prop:Wxyxy}
Let $(M^m,g)$ be a Riemannian manifold, $m\geq4$\,. The following assertions are equivalent.\\
\indent
{\rm (i)} $(M^m,g)$ is conformally-flat.\\
\indent
{\rm (ii)} $R(X,Y,X,Y)=0$ for any $X,Y\in TM$ spanning an isotropic space on $(M,g)$\,,
where $R$ is the curvature tensor of $(M,g)$\,, and $TM$ now denotes the complexified
tangent bundle.
\end{prop}
\begin{proof}
Clearly, assertion (ii) is equivalent to $W(X,Y,X,Y)=0$ for any $X,Y\in TM$ spanning an
isotropic space on $(M,g)$\,, where $W$ is the Weyl tensor of $(M,g)$\,. Therefore,
by Theorem \ref{thm:Weyl}\,, we have (i)$\Longrightarrow$(ii)\,.\\
\indent
Suppose that (ii) holds and let $(X_1,\ldots,X_m)$ be an orthonormal frame on $(M^m,g)$\,.
Then for any distinct $i,j,k,l\in\{1,\ldots,m\}$ we have
$$W(X_i\pm{\rm i}X_j,X_k+{\rm i}X_l,X_i\pm{\rm i}X_j,X_k+{\rm i}X_l)=0\;.$$
This is equivalent to the following two relations 
\begin{equation} \label{e:Wxyxy1}
W(X_i,X_k+{\rm i}X_l,X_i,X_k+{\rm i}X_l)=W(X_j,X_k+{\rm i}X_l,X_j,X_k+{\rm i}X_l)\;,
\end{equation}
\begin{equation} \label{e:Wxyxy1'}
W(X_i,X_k+{\rm i}X_l,X_j,X_k+{\rm i}X_l)=0\;.
\end{equation}
\indent
Also, by applying condition (2) of the definition of the Weyl tensor, we obtain
\begin{equation} \label{e:Wxyxy2}
\sum_{r=1}^mW(X_r,X_k+{\rm i}X_l,X_r,X_k+{\rm i}X_l)=0\;.
\end{equation}
\indent
{}From \eqref{e:Wxyxy1} and \eqref{e:Wxyxy2}\,, it follows that
$W(X_j,X_k+{\rm i}X_l,X_j,X_k+{\rm i}X_l)=0$ and, hence, $W_{jkjk}=W_{jljl}$\,,
for any distinct $j,k,l\in\{1,\ldots,m\}$\,. Therefore, for any distinct
$i,j\in\{1,\ldots,m\}$\,, we have
$$(m-1)W_{ijij}=\sum_{r=1}^mW_{irir}=0\;.$$
\indent
{}From \eqref{e:Wxyxy1'} we obtain that, for any distinct $i,j,k,l\in\{1,\ldots,m\}$\,,
we have
\begin{equation} \label{e:Wxyxy3}
\begin{split}
W_{ikjk}&=W_{iljl}\;,\\
W_{ikjl}&=-W_{iljk}\;.
\end{split}
\end{equation}
\indent
The first relation of \eqref{e:Wxyxy3} implies $W_{ikjk}=0$\,, whilst from the
second relation of \eqref{e:Wxyxy3} and the algebraic Bianchi identity
it follows quickly that $W_{ijkl}=0$\,, for any distinct $i,j,k,l\in\{1,\ldots,m\}$\,.\\
\indent
Thus, if (ii) holds then $W=0$ which, by Theorem \ref{thm:Weyl}\,, is equivalent to (i)\,.
\end{proof}

\indent
If $\H$ is a distribution on a Riemannian manifold $(M^m,g)$ we shall denote by $I^{\H}$
the integrability tensor of $\H$, which is the $\V$-valued horizontal two-form on $M^m$ defined by
$I^{\H}\!(X,Y)=-\V[X,Y]$\,, for any horizontal vector fields $X$ and $Y$, where $\V=\H^{\perp}$.\\
\indent
Next, we prove the following:

\begin{prop} \label{prop:isoI}
Let $\phi:(M^m,g)\to(N^n,h)$ be a horizontally conformal submersion between
conformally-flat Riemannian manifolds.\\
\indent
Then $g\bigl(I^{\H}\!(X,Y),I^{\H}\!(X,Y)\bigr)=0$\,, for any horizontal vectors $X$ and $Y$
spanning an isotropic space on $(M^m,g)$\,.
\end{prop}
\begin{proof}
As both the hypothesis and the conclusion are conformally-invariant, we may suppose
that $\phi:(M^m,g)\to(N^n,h)$ is a Riemannian submersion. Then the proof
follows easily from Proposition \ref{prop:Wxyxy} and the following
well-known relation of B.~O'Neill (see \cite{BaiWoo2}\,):
\begin{equation*}
R^M(X,Y,X,Y)=\phi^*(R^N)(X,Y,X,Y)-\frac34\,g(\V[X,Y],\V[X,Y]\bigr)\;,
\end{equation*}
for any horizontal vector fields $X$ and $Y$.
\end{proof}

\begin{cor} \label{cor:isoI}
Any horizontally conformal submersion, with fibres of dimension at most two,
between conformally-flat Riemannian manifolds has integrable horizontal distribution,
if the codomain has dimension at least four.
\end{cor}
\begin{proof}
Let $\phi:(M^m,g)\to(N^n,h)$ be a horizontally conformal submersion between
conformally-flat Riemannian manifolds, $m\geq n\geq4$\,.\\
\indent
Let $x\in M$ and let $E\subseteq T_xM$ be an oriented four-dimensional subspace.
{}From Proposition \ref{prop:isoI}\,, it follows that $I^{\H}_x:\Lambda^2_+E\to\V_x$
is conformal, where $\Lambda^2_+E$ is the space of self-dual bivectors on $(E,g|_E)$\,.
As $\Lambda^2_+E$ is three-dimensional, we obtain that either $I^{\H}_x=0$
or $\dim(\V_x)\geq3$\,.
\end{proof}

\indent
We end this section with an application of Corollary \ref{cor:isoI}\,.\\
\indent
An \emph{almost CR-structure}, on a manifold $M^m$, is a section $J$ of
${\rm End}(\H)$ such that $J^2=-\,{\rm Id}_{\H}$\,, where $\H$ is some distribution
on $M^m$. Obviously, $J$ is determined by its eigenbundle corresponding to $-{\rm i}$
(or ${\rm i}$)\,. Furthermore, a subbundle $\F$ of the complexified tangent bundle
of $M^m$ is the eigenbundle corresponding to $-{\rm i}$ of an almost CR-structure
on $M^m$ if and only if $\F\cap\overline{\F}=\{0\}$\,.\\
\indent
Let $J$ be an almost CR-structure on $M^m$ and let $\F$ be its eigenbundle
corresponding to $-{\rm i}$\,; $J$ is called \emph{integrable} if for any $X,\,Y\in\G(\F)$
we have $[X,Y]\in\G(\F)$\,. A \emph{CR-structure} is an integrable almost CR-structure
(see \cite{LouPan-II}\,). If $M^m$ is endowed with a Riemannian metric $g$ then
$\F$ is isotropic, with respect to $g$\,, if and only if $J$ is orthogonal,
with respect to (the Riemannian metric induced on $\H$ by) $g$\,.\\
\indent
For example, any oriented two-dimensional distribution $\V$, on a Riemannian manifold
$(M^m,g)$\,, determines two orthogonal CR-structures on $(M^m,g)$\,; at each point
$x\in M$, these are given by the rotations of angles $\pm\,\p/2$ on $\V_x$
(cf.\ \cite{Woo-4d}\,).\\
\indent
Let $\phi:(M^{n+2},g)\to(N^n,h)$ and $\psi:(N^n,h)\to(P^2,k)$ be
horizontally conformal submersions, $n\geq2$\,. Let $\V={\rm ker}\dif\!\phi$\,,
$\H=\V^{\perp}$ and let $\mathscr{K}\subseteq\H$ be the horizontal lift of
$({\rm ker}\dif\!\psi)^{\perp}$\,. Assume $\V$ and $P^2$ oriented and orient
$\mathscr{K}$ such that the isomorphism $\mathscr{K}=(\psi\circ\phi)^*(TP)$
to be orientation preserving.\\
\indent
Then the positive/negative orthogonal CR-structures determined by $\V$ and the
positive orthogonal CR-structure determined by $\mathscr{K}$ sum up to give orthogonal
almost CR-structures $J_{\pm}^{\phi,\psi}$ on $(M^{n+2},g)$\,. Obviously,
if we endow $(P^2,k)$ with its positive Hermitian structure $J^P$ then
$\psi\circ\phi:(M^{n+2},J_{\pm}^{\phi,\psi})\to(P^2,J^P)$ is holomorphic;
that is, the differential of $\psi\circ\phi$ intertwines $J_{\pm}^{\phi,\psi}$
and $J^P$.\\
\indent
We call $J_{\pm}^{\phi,\psi}$ the \emph{positive/negative almost CR-structures associated
to $\phi$ and $\psi$}\,.

\begin{prop} \label{prop:n-Wood}
Let $\phi:(M^{n+2},g)\to(N^n,h)$ be an $n$-harmonic morphism from
a Riemannian manifold of constant curvature to a conformally-flat Riemannian manifold,
and let $\psi:(N^n,h)\to(P^2,k)$ be a horizontally conformal submersion
which is real-analytic in flat local coordinates, $(n\geq4)$\,; assume
$\V\,(={\rm ker}\dif\!\phi)$ and $P^2$ oriented.
Denote by $J_{\pm}^{\phi,\psi}$ the almost CR-structures associated to $\phi$ and $\psi$\,.\\
\indent
Then either $J_+^{\phi,\psi}$ or $J_-^{\phi,\psi}$ is integrable and parallel along the
fibres of $\phi$\,.
\end{prop}
\begin{proof}
As the $n$-Laplacian on $n$-dimensional riemannian manifolds is conformally invariant,
we may suppose $(N^n,h)$ real-analytic, in flat local coordinates. Therefore, also,
$\phi$ is real-analytic.\\
\indent
Note that $\phi$ has minimal fibres \cite{Lou-p}\,. Also, by
Corollary \ref{cor:isoI}\,, the distribution $\H$ is integrable.\\
\indent
Let $\F_{\pm}$ be the eigenbundles of $J_{\pm}^{\phi,\psi}$ corresponding to $-{\rm i}$\,.
Let $Y$ be a basic vector field which locally generates $\F_+\cap\F_-$\,.
Then $Y$ is isotropic. Moreover, from the fact that $\phi$ and $\psi$ are
horizontally conformal, it follows that $\nabla_YY$ is proportional to $Y$,
where $\nabla$ is the Levi-Civita connection of $(M^{n+2},g)$\,.\\
\indent
There exists an isotropic vertical vector field $U$ such that $\F_+$ and $\F_-$
are, locally, generated by $\{U,Y\}$ and $\{\overline{U},Y\}$\,, respectively;
we may suppose that $g(U,\overline{U})=1$\,. As $Y$ is basic, $[U,Y]$ and $[\overline{U},Y]$
are vertical. Thus, $\F_+$ and $\F_-$ are integrable if and only if $g([U,Y],U)=0$
and $g([\overline{U},Y],\overline{U})=0$\,, respectively.\\
\indent
As $(M^{n+2},g)$ has constant curvature, $R^M(U,Y,Y,\overline{U})=0$\,. On the other hand,
a straightforward calculation shows that (cf.\ \cite{LouPan}\,)
\begin{equation} \label{e:n-Wood}
R^M(U,Y,Y,\overline{U})=g([U,Y],U)g([\overline{U},Y],\overline{U})\;.
\end{equation}
\indent
The proof follows.
\end{proof}

\begin{rem}
1) If $n=2$ then the conclusion of Proposition \ref{prop:n-Wood} holds
under the assumption that $(M^4,g)$ is Einstein \cite{Woo-4d}
(see \cite{LouPan} for a generalization of this result to Einstein--Weyl spaces).\\
\indent
2) Proposition \ref{prop:n-Wood}\,, also, holds under the assumption that
$\phi$ is a real-analytic horizontally conformal
submersions such that the mean curvature of $\V$ takes values
in $(\V\oplus\mathscr{K})^{\perp}$\,. Also, note that, in the proof, we have not use
the fact that $\dif\!\phi(\mathscr{K})^{\perp}\,\bigl(={\rm ker}\dif\!\psi\bigr)$ is integrable.
\end{rem}

\section{The main result} \label{section:main}

\indent
This section is devoted to the following result and its consequences.

\begin{thm} \label{thm:cf1d}
Let $\phi:(M^{n+1},g)\to(N^n,h)$ be a harmonic morphism between Riemannian manifolds,
$n\geq3$\,; denote by $\l$ the dilation of $\phi$\,.\\
\indent
If $(M^{n+1},g)$ is real-analytic and conformally-flat then \emph{either}\\
\indent
\quad{\rm (i)} $\phi$ is of Killing type, \emph{or}\\
\indent
\quad{\rm (ii)} the horizontal distribution of $\phi$ is integrable and its
leaves endowed with the metrics induced by $\l^{-2n+4}g$ have constant curvature.
\end{thm}
\begin{proof}
\indent
By a result of \cite{PanWoo-d}\,, at least away of the critical points (which may
occur only if $n=3$\,, see \cite{BaiWoo2}\,), we have
$\phi:(M^{n+1},g)\to(N^n,h)$ real-analytic.\\
\indent
As the dimension of the intersection of (the complexification of) $\H$
with any isotropic two-dimensional space, on $(M^{n+1},g)$\,, is at least $1$\,,
Proposition \ref{prop:Wxyxy} implies that $(M^{n+1},g)$ is
conformally-flat if and only if, for any $U\in\G(\V)$ and $X,Y\in\G(\H)$
with $g(U,U)=g(X,X)$\,, $g(X,Y)=0$\,, $g(Y,Y)=0$\,, we have
$R^M(U\pm{\rm i}X,Y,U\pm{\rm i}X,Y)=0$\,; equivalently,
\begin{equation} \label{e:cf1d1}
\begin{split}
R^M(U,Y,U,Y)&=R^M(X,Y,X,Y)\\
R^M(U,Y,X,Y)&=0\;.
\end{split}
\end{equation}
{}From \eqref{e:xyzv}\,, it follows quickly that the second relation of \eqref{e:cf1d1}
is equivalent to
\begin{equation} \label{e:cf1d1second}
(^h\nabla_Y\O)(X,Y)+3(n-1)Y(\s)\O(X,Y)=0\;.
\end{equation}
Thus, by assuming $X$ and $Y$ basic and using Proposition \ref{prop:Bryant}\,,
we obtain
\begin{equation} \label{e:cf1d2}
Y(V(\s))\O(X,Y)=0\;,
\end{equation}
where $V$ is the fundamental vector field of $\phi$\,.\\
\indent
Next, we shall use the first relation of \eqref{e:cf1d1}\,. For this, we assume
$X$ and $Y$ basic with $g(X,X)=e^{-2\s}$ (equivalently, $h(X,X)=1$\,), and
$U=e^{-(n-1)\s}V$ (so that, $g(U,U)=g(X,X)$\,). Thus, the first relation
of \eqref{e:cf1d1} becomes $$e^{-(2n-2)\s}R^M(V,Y,V,Y)=R^M(X,Y,X,Y)$$
which, by applying \eqref{e:xvyv} and \eqref{e:xyzh}\,, is equivalent to
\begin{equation} \label{e:cf1d3}
\begin{split}
R^N(X,Y,X,Y)=-(n-1)h(^h&\nabla_Y(\H(\grad_h\s)),Y)-(n-1)^2Y(\s)^2\\
&+\frac14\,e^{(2n-2)\s}\bigl\{h(i_Y\O,i_Y\O)+3\O(X,Y)^2\bigr\}\;,
\end{split}
\end{equation}
where we have denoted by the samy symbol $R^N$ and its pull-back by $\phi$ to $M^{n+1}$.\\
\indent
We may assume that $Y$ is the horizontal lift of an isotropic geodesic (local)
vector field on (the complexification of) $(N^n,h)$\,; equivalently, $^h\nabla_YY=0$\,.
Then \eqref{e:cf1d3} becomes
\begin{equation} \label{e:cf1d3'}
\begin{split}
R^N(X,Y,X,Y)=-(n-1)&Y(Y(\s))-(n-1)^2Y(\s)^2\\
&+\frac14\,e^{(2n-2)\s}\bigl\{h(i_Y\O,i_Y\O)+3\O(X,Y)^2\bigr\}\;.
\end{split}
\end{equation}
\indent
As $R^N(X,Y,X,Y)$ is basic, from \eqref{e:cf1d2} and \eqref{e:cf1d3'} it easily follows
that either $\O=0$ or
\begin{equation} \label{e:cf1d4}
V(\s)\bigl\{h(i_Y\O,i_Y\O)+3\O(X,Y)^2\bigr\}=0\;.
\end{equation}
\indent
Now, from $\O\neq0$ it follows that there exist $Y\in\H$ isotropic
and $X\in Y^{\perp}\cap\H$ such that the second factor of the left hand side of
\eqref{e:cf1d4} is not zero. Thus, we have proved that either $\O=0$ (equivalently,
$\H$ is integrable) or $V(\s)=0$ (equivalently, $\phi$ is of Killing type).\\
\indent
Next, we study the case $\O=0$\,. Then \eqref{e:cf1d1second} (and hence, also, the
second relation of \eqref{e:cf1d1}\,) is automatically satisfied, whilst \eqref{e:cf1d3}
is equivalent to
\begin{equation} \label{e:cf1d5}
R^N(X,Y,X,Y)=\,^h\nabla(\dif^{\H}\!u)(Y,Y)-(\dif^{\H}\!u)(Y)^2\;,
\end{equation}
where $u=-(n-1)\s$ and, recall that, $X$ and $Y$ are basic with $h(X,X)=1$\,, $h(X,Y)=0$
and $h(Y,Y)=0$\,.\\
\indent
Let $h_1=e^{2u}h|_{\H}=e^{(-2n+4)\s}g|_{\H}$\,.\\
\indent
We have proved that, if $\H$ is integrable, \eqref{e:cf1d1} is equivalent
to the fact that the curvature tensor $R^P$ of any leaf $P$ of $\H$\,, endowed with the
metric induced by $h_1$\,, satisfies $R^P(X,Y,X,Y)=0$\,.\\
\indent
It follows that if $\H$ is integrable then $h_1$ induces a conformally-flat Einstein metric
on each leaf of $\H$\,; equivalently, $h_1$ induces a metric of constant curvature
on each leaf of $\H$\,. The proof is complete.
\end{proof}

\begin{exm}
Let $(N^n,h)$ be $\R^n$, endowed with the canonical metric, and let
$$M^{n+1}=\bigl\{(t,x)\in\R\times\R^n\,\big{|}\,|tx|<1\bigr\}\;,$$
where $|\cdot|$ denotes the Euclidean norm on $\R^n$.\\
\indent
Define $\l:M^{n+1}\to(0,\infty)$ by $\l(t,x)=(1-|tx|^2)^{\frac{1}{n-1}}$, $(t,x)\in M^{n+1}$,
and let $g=\l^{-2}h+\l^{2n-4}\dif\!t^2$.\\
\indent
Then $\phi:(M^{n+1},g)\to(N^n,h)$\,, $(t,x)\mapsto x$\,, is a harmonic morphism which
satisfies assertion (ii) of Theorem \ref{thm:cf1d}\,; in particular, $(M^{n+1},g)$
is conformally-flat, $(n\geq3)$\,. Furthermore,
$\phi$ is neither of Killing type nor its fibres are geodesics.
\end{exm}

\begin{rem}
If $n=3$ then Theorem \ref{thm:cf1d} holds, also, in the complex-analytic category.
Indeed, the only point in the proof
of Theorem \ref{thm:cf1d} where it is essential for $\phi$
to be `real' is when we deduce from $\O\neq0$ that there exist $Y\in\H$
isotropic and $X\in Y^{\perp}\cap\H$ such that the second factor of the left
hand side of \eqref{e:cf1d4} is not zero. But, if $n=3$ and $h(X,X)=1$ then
$$h(i_Y\O,i_Y\O)+3\O(X,Y)^2=4\O(X,Y)^2\;,$$
which, also, in the complex-analytic category, is not zero, for suitable
choices of $X$ and $Y$, if $\O\neq0$\,.
\end{rem}

\indent
Next, we discuss the case when both the domain and codomain, of a harmonic
morphism with one-dimensional fibres, are conformally-flat; the notations
are as in Section 1.

\begin{cor} \label{cor:cf1d}
Let $\phi:(M^{n+1},g)\to(N^n,h)$ be a submersive harmonic morphism
with connected one-dimensional fibres, $(n\geq3)$\,.\\
\indent
The following assertions are equivalent.\\
\indent
\quad{\rm (i)} $(M^{n+1},g)$ and $(N^n,h)$ are conformally-flat.\\
\indent
\quad{\rm (ii)} One of the following assertions holds:\\
\indent
\quad\quad{\rm (iia)} $\phi$ is of Killing type, $n=3$\,, and, up to a
homothety, $\O$ is the volume form of a Riemannian foliation by
geodesic surfaces, of sectional curvature $1$\,, on $(N^3,\l^{-4}h)$\,.\\
\indent
\quad\quad{\rm (iib)} The horizontal distribution of $\phi$ is integrable and its
leaves endowed with the metrics induced by $\l^{-2n+4}g$ have constant curvature.
\end{cor}
\begin{proof}
If $n\geq4$ this follows from Corollary \ref{cor:isoI} and the proof
of Theorem \ref{thm:cf1d}\,.\\
\indent
Assume $n=3$\,. Then by the proof of Theorem \ref{thm:cf1d}\,, if $(M^4,g)$ is conformally-flat,
on each connected component of a dense open subset of $M^4$,
either $\phi$ is of Killing type or (iib) holds.\\
\indent
If $\phi$ is of Killing type then there exist Weyl connections $D_{\pm}$
on $(N^3,[h])$ such that $\phi:(M^4,[g])\to(N^3,[h],D_{\pm})$ is $\pm$twistorial,
in the sense of \cite[Example 4.8]{LouPan-II} (cf.\ \cite{PanWoo-sd}\,). Furthermore,
$(M^4,g)$ is conformally-flat if and only if both $D_{\pm}$ are Einstein--Weyl
(see \cite{PanWoo-sd} and the references therein).
Also, if $(N^3,h)$ is conformally-flat then $D_{\pm}$ are Einstein--Weyl
if and only if, locally, $D_{\pm}$ are the Levi-Civita connections of
constant curvature representatives $h_{\pm}$ of $[h]$
(see \cite{Cal-F}\,).\\
\indent
We claim that if $n=3$ and $\phi$ is of Killing type then, with the same notations
as above, the following assertions are equivalent:\\
\indent
\quad(a) Up to a homothety, $\O$ is the volume form of a Riemannian foliation
by geodesic surfaces, of sectional curvature $1$\,, on $(N^3,\l^{-4}h)$\,.\\
\indent
\quad(a$'$) $D_{\pm}$ are, locally, the Levi-Civita connections of constant curvature
representatives $h_{\pm}$ of $[h]$\,, and $D_+\neq D_-$\,.\\
\indent
Indeed, if $\phi$ is of Killing type then, by replacing $g$ and $h$ with $\l^{-2}g$
and $\l^{-4}h$\,, respectively, we may suppose that $\phi$ is a Riemannian submersion
with geodesic fibres. Then the Lee forms $\a_{\pm}$ of $D_{\pm}$\,, with respect to $h$\,,
are given by $\a_{\pm}=\pm*_h\O$ (see \cite{LouPan-II}\,, \cite{PanWoo-sd}\,),
where $*_h$ is the Hodge $*$-operator of $h$\,,
with respect to some local orientation, and we have denoted by the same symbol $\O$ and
the two-form on $N^3$ whose pull-back by $\phi$ is $\O$\,. Hence, if (a$'$) holds then
$h_{\pm}=e^{\pm2u}h$ where $u$ is characterised by $\dif\!u=*_h\O$\,;
in particular, $u$ is a harmonic (local) function on $(N^3,h)$\,.\\
\indent
It follows that (a$'$) is equivalent to the following assertion:\\
\indent
\quad(a$''$) Locally, there exists a nonconstant function $u$ on $N^3$ such that
$$\dif\!u=*_h\O\,,\;(\nabla^h\!\dif\!u)(Y,Y)=0\,,\;{\rm Ric}^h(Y,Y)=-\dif\!u(Y)^2\,,$$
for any isotropic vector $Y$ on $(N^3,h)$\,, where $\nabla^h$
is the Levi-Civita connection of $(N^3,h)$ and ${\rm Ric}^h$ is the Ricci tensor
of $(N^3,h)$\,.\\
\indent
Now, by applying, for example, Lemma \ref{lem:curv}\,, we obtain that assertion (a$''$)
is equivalent to the following:\\
\indent
\quad(a$'''$) Locally, there exists a nonconstant function $u$ on $N^3$ such that
$$\dif\!u=*_h\O\,,\;\nabla^h\!\dif\!u=0$$
and the level surfaces of $u$ have sectional curvature equal to $|\!\dif\!u|^2$.\\
\indent
The proof of (a)$\iff$(a$'$) follows.\\
\indent
We have thus proved that (ii)$\Longrightarrow$(i)\,, and if (i) holds then, also, (ii) holds
on each connected component of a dense open subset
of $M^4$.\\
\indent
To complete the prof of (i)$\Longrightarrow$(ii)\,, define a connection $\nabla$ on $\H$ by
$$\nabla_EX=\H\nabla^{\l^{-2}g}_{\H E}X+\H[\V E,X]$$
for any vector field $E$ and horizontal vector field $X$\,, where
$\nabla^{\l^{-2}g}$ is the Levi-Civita connection of $(M^4,\l^{-2}g)$\,.\\
\indent
If we assume (i) then, from the fact that (ii) holds on each connected component of a dense open subset
of $M^4$, it follows quickly that $\nabla\O=0$\,. Therefore either $\O$ is nowhere zero or $\O=0$ on $M^4$.
The proof is complete.
\end{proof}

\begin{exm} \label{exm:Hopfmap}
Let $\p:\R^4\to\R^3$ be the \emph{Hopf polynomial map} defined by
$\p(z_1,z_2)=(|z_1|^2-|z_2|^2,2z_1\overline{z_2})$\,, for any $(z_1,z_2)\in\R^4\,(=\C^{\!2})$\,.\\
\indent
Then $\p|_{\R^4\setminus\{0\}}$ satisfies assertion (iia) of Corollary \ref{cor:cf1d}\,.
\end{exm}

\indent
We end with the following consequence of Corollary \ref{cor:cf1d}\,.

\begin{cor} \label{cor:cf1dHopf}
The Hopf polynomial map $\p:\R^4\to\R^3$ is, up to local conformal diffeomorphisms with basic conformality factors,
the only harmonic morphism with one-dimensional fibres and nonintegrable horizontal distribution between
conformally-flat Riemannian manifolds, of dimensions at least three.
\end{cor}
\begin{proof}
This follows from the fact that any harmonic morphism which satisfies assertion (iia) of
Corollary \ref{cor:cf1d} is, locally, the Hopf polynomial map $\p:\R^4\to\R^3$, up to conformal
changes with basic factor.
\end{proof}

\end{document}